\magnification=\magstep1

\def\pne{{\cal P}({}^nE)}
\def\dis{\displaystyle}
\def\o{\over}
\def\ons#1{\widehat{\bigotimes
_{n,s,#1}}}
\null\vskip1.5truein

\centerline{{\bf A Dvoretsky Theorem for Polynomials}}
\centerline{{\bf by}}
\centerline{{\bf Se\'an Dineen\footnote{*}{\sevenrm
The author  thanks Universidade
Federal do  Rio de Janeiro for support during the period
this research was initiated.} (University College Dublin)}}

\bigskip\bigskip
\centerline{\bf Abstract}
{\narrower\narrower\narrower
\noindent  We lift upper and lower estimates from
linear functionals to $n$-homogeneous polynomials and using
this result show that $l_\infty$ is finitely represented in
the space of $n$-homogeneous polynomials, $n\ge2$, for any
infinite dimensional Banach space.
Refinements are also given.\bigskip}

\noindent The classical Dvoretsky spherical sections theorem
[5,13] states that $l_2$ is finitely represented in any
infinite dimensional Banach space. Using this, the Riesz
Representation theorem (for finite dimensional  $l_p$
spaces) and the Hahn-Banach theorem we show that $l_\infty$
is finitely represented in ${\cal P}({}^nE)$, for any
infinite dimensional Banach space and any $n\ge2$.  This
shows that $\pne$ does not have any non-trivial
superproperties and explains why spaces such as Tsirelson's
space play such a positive role in the recent theory of
polynomials on Banach spaces
$\bigl($[1,2,6,7,8,9,10]$\bigr)$.  We refer to [3,11,12] for
properties of Banach spaces and to [4] for properties of
polynomials.

\proclaim Theorem 1.  Suppose $E$ is a
Banach space, $1<p\le\infty$, $\{\phi_j\}_{j=1}^k$ is a finite sequence of
vectors in $E'$ and $A$ and $B$ are positive constants such
that
$$A^p\sum_{j=1}^k|\alpha_j|^p
\le\left\|\sum_{j=1}^k\alpha_j\phi_j\right\|^p \le
B^p\sum_{j=1}^k|\alpha_j|^p\eqno(1)$$
for any sequence of scalars
$(\alpha_j)_{j=1}^k$.  Then for any integer $n$, $n\ge q$
where $\dis{1\o p}+{1\o q}=1$, and any sequence of scalars
$(\alpha_j)_{j=1}^k$ we have
$$A^n\sup_{1\le j\le
k}|\alpha_j|\le\left\|\sum_{j=1}^k\alpha_j\phi{}_j^n\right\|
\le B^n\sup_{1\le j\le k}|\alpha_j|\eqno(2)$$

\noindent{\bf Proof}.  For any $x\in E$, $\|x\|\le1$, we have
$$\sup_{\sum_{j=1}^k|\alpha_j|^p\le1}
\left|\sum_{j=1}^k\alpha_j\phi_j(x)\right|^p\le B^p.$$
Since $(l_p)'=l_q$ this implies
$$\sup_{\|x\|\le1}\sum_{j=1}^k\bigl|\phi_j(x)\bigr|^q\le
B^q.$$
If $n\ge q$ then
$$\eqalignno{
\sup_{\|x\|\le1}\left|\sum_{j=1}^k\alpha_j\phi_j^n(x)\right|
&\le\sup_j|\alpha_j|\cdot
B^n\sup_{\|x\|\le1}\sum_{j=1}^k\left|{\phi_j(x)\o
B}\right|^n\cr
&\le B^n\sup_j|\alpha_j|.&(3)\cr}$$
On the other hand
$$A^p\sum_{j=1}^k|\alpha_j|^p=A^p
\sup_{\sum_{j=1}^k|\beta_j|^q\le1}
\left|\sum_{j=1}^k\alpha_j\beta_j\right|^p
\le\sup_{\|x\|\le1}
\left|\sum_{j=1}^k\alpha_j\phi_j(x)\right|^p$$
Since the set
$\displaystyle\Bigl\{\bigl(\phi_j(x)\bigr)_{j=1}^k;
\|x\|\le1\Bigr\}$
is a convex balanced set, the Hahn-Banach theorem
implies that
$$A\cdot
B_{l_q^k}=A\cdot\bigg\{(\beta_j)_{j=1}^k;
\sum_{j=1}^k|\beta_j|^q\le1\biggr\}
\subset\overline{\biggl\{\bigl(\phi_j(x)\bigr)_{j=1}^k;
\|x\|\le1\biggr\}}$$
Hence, for any fixed integer $l$, $1\le l\le k$, there
exists $(x_n)_n$ in $E$, $\|x_n\|\le1$, such that
$$\phi_l(x_n)\to A\hbox{ and, for }j\ne l,\quad
\phi_j(x_n)\to 0\hbox{ as }n\to\infty.$$
This implies that
$$\sup_{\|x\|\le1}\left|\sum_{j=1}^k\alpha_j\phi_j^n(x)\right|
\ge A^n\cdot\sup_{1\le j\le k}|\alpha_j|\eqno(4)$$
and the inequalities (3) and (4) prove the proposition.
\bigskip
Note that in the proof of theorem~1 we have actually shown
that the left (resp.~right) hand side of the inequality~(1)
implies the left (resp.~right) hand side of the
inequality~(2).  Conditions (1) and (2) can be rephrased in
terms of the Banach-Mazur distance $d$ to give the following
result.
\proclaim Corollary 2.  Let $E$ denote a $k$-dimensional
Banach space and suppose $d(E,l_p^k)\le C$ where $1\le
p<\infty$.  Then, for $n\ge p$, $\pne$ contains a
$k$-dimensional subspace $F$ such that
$$d(F,l_\infty^k)\le C^n.$$

\proclaim Corollary 3.  If $E$ is an infinite dimensional
Banach space and $n\ge2$ then $l_\infty$ is finitely
represented in $\pne$.

\noindent{\bf Proof}.  By the classical Dvoretsky
theorem we can choose for any positive integer $k$ and any
$\epsilon>0$ vectors $\{\phi_j\}_{j=1}^k$ in $E'$ such that
$$\sum_{j=1}^k|\alpha_j|^2
\le\left\|\sum_{j=1}^k\alpha_j\phi_j\right\|^2
\le(1+\epsilon)^2\sum_{j=1}^k|\alpha_j|^2$$
for any sequence of scalars $(\alpha_j)_{j=1}^k$.  Hence, by
theorem~1, we have, for $n\ge2$ and any
$(\alpha_j)_{j=1}^k$,
$$\sup_{1\le j\le k}|\alpha_j|\le\left\|\sum_{j=1}^k
\alpha_j\phi_j^n\right\|\le(1+\epsilon)^n\sup_{1\le
j\le k}|\alpha_j|$$
This proves the corollary.
\bigskip
In fact it is easily seen that the above shows that
$l_\infty$ is finitely represented in the space of
polynomials of finite type.

\proclaim Corollary 4.  If $l_p$, $1\le p<\infty$, is
a quotient of $E$ then
$l_\infty$ is a subspace of $\pne,n\ge p$, and $l_1$ is a
complemented subspace of the completed symmetric tensor
product endowed with the projective topology, $\displaystyle
\ons \pi E$.

\noindent{\bf Proof}.  Let $\dis{1\o p}+{1\o q}=1$.  We can
choose constants $A$ and $B$ independent of $k$ and vectors
$(\phi_k)_k$ in $E'$ such that
$$A^q\cdot\sum_{j=1}^k|\alpha_j|^q
\le\left\|\sum_{j=1}^k\alpha_j\phi_j\right\|^q
\le B^q\cdot\sum_{j=1}^k|\alpha_j|^q$$
for any sequence of scalars $(\alpha_j)_j$.  Theorem~1
implies that, for $n\ge p$,
$$A^n\sup_{1\le j\le k}|\alpha_j|
\le\left\|\sum_{j=1}^k\alpha_j\phi_j^n\right\|
\le B^n\sup_{1\le j\le k}|\alpha_j|$$
for any sequence of scalars $(\alpha_j)_j$ and any $k$.
Hence $\{\phi_j^n\}_{j=1}^\infty$ is equivalent to the unit
vector basis of $c_0$.  Since $\pne$ is a dual space this
implies $l_\infty\hookrightarrow\pne$.  Since
$\displaystyle\left(\ons \pi E\right)'\cong\pne$
it follows by
[3, p.~48] that $l_1$ is a complemented subspace of
$\displaystyle
{\ons \pi}E$.  This completes the proof.
\bigskip
This result for $p=2$
is given in [8, proposition~13]
and also implies the well known fact that
${\cal P}({}^nl_p)$ is
not reflexive if $n\ge p$.

We now extend the result given in corollary~3 and at the
same time obtain a refinement of [6, theorem~1(ii)].
The elements of $\displaystyle\ons \epsilon E'$ are
$n$-homogeneous polynomials on $E$ which are uniformly
weakly continuous on bounded subsets of~$E$.  Hence they
have unique extensions to $E''$.  We use the notation
$\tilde P$ to denote this extension.

\proclaim Lemma 5.  A bounded sequence $(P_j)_j$ in
$\displaystyle\ons
\epsilon E'$ is a weakly null sequence if and only if
$\tilde P_j(x'')\to0$ as $j\to\infty$ for any $x''\in E''$.

\noindent{\bf Proof}.  If $\phi\in\left(
\displaystyle\ons \epsilon E'\right)'$ then there exists a
regular Borel measure $\mu$ on
$\left(\overline{B}_{E''},\sigma(E'',E')\right)$ such that
$$\phi(P)=\int_{\overline{B}_{E''}}\tilde
P(x'')\,d\mu(x'')$$
for all $\displaystyle P\in\ons \epsilon E'$.

If $(P_j)_j$ is bounded then $(\tilde P_j)_j$ is uniformly
bounded on $\overline{B}_{E''}$.  If $\tilde P_j(x'')\to 0$
as $j\to\infty$ for each $x''$ in $\overline{B}_{E''}$ then
the Lebesgue dominated convergence theorem implies that
$\phi(P_j)\to0$ as $j\to\infty$.  Hence $(P_j)_j$ is weakly
null.  The converse is obvious.

\proclaim Proposition 6.  If $E$ is a Banach space and $E'$
contains a weakly null sequence of unit vectors, which
satisfies an upper $q$ estimate, $q<\infty$, then for $n\ge
p$, $\displaystyle{1\o p}+{1\o q}=1$, we have
$l_\infty\hookrightarrow\pne$.

\noindent{\bf Proof}.  If suffices to show
$\displaystyle c_0
\hookrightarrow{\ons \epsilon}E'$.
Let $(\phi_j)_j$ denote a weakly null sequence of unit
vectors in $E'$ which satisfies an upper $q$-estimate. Then
$(\phi_j^n)_j$ is a sequence of unit vectors in
$\displaystyle\ons \pi E'$. Since $(\phi_j)_j$ is weakly
null lemma~5 implies that $(\phi_j^n)_j$ is a weakly null
sequence in $\displaystyle\ons \epsilon E'$. By the
Bessaga-Pelczynski selection principle [3, p.~42 and 11,
p.~5] the
sequence $(\phi_j^n)_j$ contains a subsequence which forms a
basic sequence.  Since upper $q$ estimates are inherited by
subsequences we may suppose that $(\phi_j^n)_j$ is a basic
sequence.

Hence there exists $B>0$ such that
$$\left\|\sum_{j=1}^k\alpha_j\phi_j\right\|^q\le
B^q\sum_{j=1}^k|\alpha_j|^q$$ for any integer $k$ and any
sequence of scalars $(\alpha_j)_j$.

By theorem~1 we have
$$\left\|\sum_{j=1}^k\alpha_j\phi_j^n\right\|\le B^n
\sup_{1\le j\le n}|\alpha_j|$$ for any sequence of scalars
$(\alpha_j)_j$.

Since $(\phi_j^n)_j$ is a basic sequence the closed subspace
of $\displaystyle\ons \epsilon E'$ spanned by
$(\phi_j^n)_j$ is isomorphic to $c_0$.  This completes the
proof.
\bigskip
Proposition~5 applies in particular to
reflexive Banach lattices which satisfy a lower $p$-estimate
$\displaystyle{1\o p}+{1\o q}=1$ ([12]).

\bigskip
\centerline{\bf Bibliography}
\item{[1]} R. Alencar, R. M. Aron, S. Dineen, {\sl A
reflexive space of holomorphic functions in infinitely many
variables},  Proc. A. M. S., {\bf90}, 3, 1984, 407--411.
\item{[2]} R. M. Aron, S. Dineen, {\sl$Q$-reflexive Banach
spaces,} preprint.
\item{[3]} J. Diestel, {\sl Sequences and Series in
Banach Spaces}, Springer, Graduate Text in Mathematics,
1984.
\item{[4]} S. Dineen, {\sl Complex Analysis in locally
convex spaces}, North Holland Math Studies, {\bf57}, 1981.
\item{[5]}  A. Dvoretzty, {\sl Some results on convex bodies
and Banach Spaces},  Proc.~International Symposium on Linear
spaces, Jerusalem, 1961, 123--160.
\item{[6]} J. Farmer, {\sl Polynomial Reflexivity in Banach
Spaces}, preprint.
\item{[7]} J. Farmer, W. B. Johnson, {\sl Polynomial Schur
and Polynomial Dunford-Pettis properties,} Contemporary
Mathematics, {\bf144}, AMS, 1993, 95--105.
\item{[8]} M. Gonzalez, J. Gutierrez, {\sl Unconditionally
converging Polynomials on Banach spaces,} preprint.
\item{[9]} R. Gonzalo, J. A. Jarmillo, {\sl Compact
Polynomials between Banach spaces}, preprint.
\item{[10]} J. Jaramillo, A. Prieto, {\sl Weak polynomial
convergence on a Banach space}, preprint.
\item{[11]} J. Lindenstrauss, L. Tzafriri, {\sl Classical
Banach spaces I}, ``Sequence Spaces'', Springer Verlag,
1977.
\item{[12]} L. Lindenstrauss, L. Tsafriri, {\sl Classical
Banach Spaces II}, Springer, Graduate Text in Mathematics,
1979.
\item{[13]}  G. Pisier, {\sl The volume of convex bodies and
Banach  space geometry}, Cambridge Tracts in Mathematics,
{\bf94}, 1989.
\bigskip
\leftline{Department of Mathematics,}
\leftline{University College Dublin,}
\leftline{Belfield,}
\leftline{Dublin 4,}
\leftline{Ireland.}
\leftline{e-mail $<$sdineen@irlearn.bitnet$>$}
\bye